\documentclass{llncs}

\newtheorem{conjectureBC}{Brouwer's Conjecture}

\newtheorem{Co}{Conjecture}
\newcommand{\mathsym}[1]{{}}
\newcommand{\unicode}[1]{{}}

\begin{document}

	\title{Proof of  Brouwer's  Conjecture (BC) For all Graphs With Number of Vertices $n>n_0$ Assuming that BC Holds for $n\leq n_0$ for some $n_0 \leq 10^{24}$} 
	\author{Vladimir Blinovsky\inst{1,2}\and Llohann D Speranca\inst{1}\and Alexander Pchelintsev\inst{3}}
	\institute{
		Universidade Federal de Sao Paulo (UNIFESP), \\
		Campus Sao Jose dos Campos. Instituto de Ciencia e Tecnologia (ICT),
		Brazil\inst{1}, \\
		Institute for Information Transmission Problems,\\
		B. Karetnyi 19, Moscow, Russia\inst{2},\\
		Department of Higher Mathematics, \\
		Tambov State Technical University, \\
		Sovetskaja Str. 106, 392000 Tambov, Russia\\
		\email{vblinovs@yandex.ru\inst{2},\\
		lsperanca@gmail.com\inst{1}, \\
		pchelintsev.an@yandex.ru\inst{3}}}
	\maketitle
	
	\bigskip

	\begin{abstract}
				
		Assuming that Brouwer's Conjecture - the upper bound for the sum of $t\leq n-1$ largest eigenvalues of Laplacian graph on $n$ vertices is true for $n\leq 10^{24}$ we prove the Brouwer`s Conjecture ($\hbox{BC}$) for $n>10^{24}$
	\end{abstract}
	\bigskip 
	\bigskip
	
	Let $A$ be $n\times n$ incidence matrix of simple undirected graph $G$:
	\begin{eqnarray*}
		a_{i,j}=\left\{\begin{array}{ll}
			1,& \hbox{iff}\  (i,j)\in G,\\
			0,& \hbox{otherwise}.
		\end{array}
		\right.
	\end{eqnarray*}
	Define the Laplacian $L(G)$ of $G$ as follows
	$$
	L(G)=D-A,
	$$
	where diagonal $n\times n$ matrix $D$ has entries
	$$
	d_{i}=|\{ j:\ (i,j)\in G\}|.
	$$
	We have $\sum_i d_{i}=2m$, were $m$ is number of edges in $G$.
	Considering $G$ as directed graph with some choice of ordering of vertices in $G$ define $m\times n$ matrix
	$B$:
	$$
	b_{i,j}=\left\{\begin{array}{ll}
		1,& \hbox{if}\ j\ \hbox{heard vertex in edge}\ $i$,\\
		-1, & \hbox{if}\ j\ \hbox{tail vertex in edge}\ $i$,\\
		0, & \hbox{otherwise}
	\end{array}
	\right.
	$$
	Then $L(G)=B^T B$ and hence eigenvalues of matrix $L(G)$ are nonnegative:
	$$
	0=\mu_n (L(G))\leq \mu_{n-1} (L(G))\leq\ldots\leq\mu_1 (L(G)).
	$$
	%	\bigskip
	
	\begin{conjectureBC}
		\cite{sa}.
		For every graph $G\subset {[n]\choose 2}$ and integer $t\in  [n-1]$, the following inequality is valid:
		$$
		S_t(G)=	\sum_{i=1}^t \mu_i (L(G))\leq m+{t+1\choose 2},\ t\in [n].
		$$		
	\end{conjectureBC}
	In this article we prove the validness of this conjecture under the assumption than it is true for all $n\leq n_0$ where $n_0 \leq 10^{24}$. 
	
	For convenience, we denote
	\[\Delta_t(G)=S_t(G)-m(G)-{t+1\choose 2}. \]
	Whenever $\Delta_t(G)\leq 0$, we say that $G$ satisfy $BC_t$.  
	
	It is known to be valid for trees~\cite{1}, for $k=1,2,n-1,n$, for unicyclic and bicyclic graphs~\cite{9}, for regular graphs~\cite{f}, for $n\leq 10$ it was checked by A. Brouwer using a computer. In~\cite{sw} was proved that Brouwer's  conjecture holds asymptotically almost surely.

	Before the proof of Brouwer`s conjecture under above assumption  (we call it below "Conjecture") we introduce some consequences of its validity.

	Define set of conjugate degrees
	$$
	d^* (G)=\{ d^*_1 ,\ldots , d^*_n\},\  d^*_ i =|\{ j: d_j \geq i\}|.
	$$
	We say that the set $E$ of edges is compressed if from $e=(i<j)\in E $ it follows that $e=(i_1 <j_1 )\in E$, were $i_1 \leq i,\ j_1 \leq j$.
	V.Ch\'{a}tal and P.Hammer in~\cite{lk} introduce the notion of  threshold  graph. It can be defined as a graph isomorphic (up to permutations on vertices $[n]$) to a graph with compressed edge set.

	Grone - Merris Conjecture~\cite{tr}, which was proved by Bai~\cite{m}, we call it GMB, theorem says that the following upper bound is valid
	\begin{equation}
	\label{d1}
	\sum_{i=1}^t \mu_i (L(G))\leq \sum_{i=1}^t d^*_i (G).
	\end{equation}
	It is known~\cite{2} that for threshold graphs there is equality in the last relation.

	Say that a graph $G$ on $n$ nodes with $m=m(G)$ edges is spectrally threshold dominated~(\cite{4}) if for each $t\in [n]$  there is a threshold graph $\hat{G}$ having the same number of nodes and edges satisfying $$\sum_{i=1}^t \mu_i (L(G))\leq\sum_{i=1}^t \mu_i (L(\hat{G}))=\sum_{i=1}^t d^*_i (L(\hat{G})).$$
	In paper~\cite{4} Helmberg and Trevisan  proved
	the following
	\begin{Co}
	
	 Graph $G$ is  spectrally threshold dominated iff Conjecture for this graph is valid.
	 \label{co22}
	 \end{Co}
	  We introduce here their  proof via construction of the set of conjugate degrees of optimal threshold graphs. 
	
	We construct for arbitrary $n,m=m(G), t$ threshold graph $T$ that attains Brouwer's bound for the sum of eigenvalues. Denote by $\hbox{Tr}(n,m)$ the set of threshold graphs with $n$ vertices and $m$ edges. To each graph with degree sequence $d_i \geq d_{i+1}$ define Ferrers diagram of $n$ rows, s.t.  $i-$th row displays
$d_i$ boxes aligned to the left. 	          

Next we demonstrate for arbitrary $t\in [n]$  that $\min\{tn, m(G) + t(t +
1)/2, 2m(G)\} = \max_{T\in \hbox{Tr}(n,m)}   \sum_{i=1}^t d^*_i (T)$.  This together with~(\ref{d1}) deliver the proof of Conjecture~\ref{co22}. 

Depending on
the relation between $t, n$ and $m(G)$, we consider the following cases:
\begin{enumerate}
\item
Case 1. $\min\{ tn, m(G) + t(t + 1)/2, 2m(G)\} = tn$. Consider the threshold graph $T$ constructed
by filling up the Ferrers diagram below the diagonal in column wise order (on and above the

diagonal in corresponding row wise order). The first $t$ columns below the diagonal are fully
filled because they require $tn- t(t + 1)/2 \leq  m(G)$ boxes. Hence $T$ satisfies $d^*_i (T)=n$
for $i\in [t]$  and $\sum_{i=1}^t d^* _i (T)=tn$. This is the maximum attainable over all
threshold graphs on $n$ nodes.
\item

Case 2. $\min\{tn, m(G) + t(t + 1)/2, 2m(G)\} = m(G) + t(t + 1)/2$: In this case put $h = \lfloor \frac{m(G)}{t}+\frac{t+1}{2}\rfloor<n$ and $r=m(G)+t(t+1)-th< t$. 
Note that this implies $h \geq t + 1$. Define a
threshold graph $T$ on $n$ nodes with $m(T)=m(G)$ edges of trace $t$ by the conjugate degrees
$$
d^*_i (T)=\left\{\begin{array}{ll} h+1,\ & i\leq r,\\
h, & r<i\leq t,
\end{array}
\right.
$$

then $\sum_{i=1}^t \lambda_i (T) = \sum_{i=1}^t d^*_i (T)=m(T) + t(t + 1)/2$. This value cannot be exceeded by any
threshold graph on $n$ nodes with $m$ edges by the GMB  theorem,
because in the Ferrers diagram of the conjugate degrees up to column $t$ all boxes are used on
and above the diagonal, while all possible m boxes are included below the diagonal.

\item

Case 3. $\min\{tn, m(G) + t(t + 1)/2, 2m(G)\} = 2m(G)$. Put $h= \max\{h\in [n] : h(h + 1) \leq
2m(G)\} < t$ and $r = (2m(G) -h(h + 1))/2 < h + 1$, then the threshold graph $T$ of trace $h$ with
conjugate degrees
$$
d^*_i(T) =\left\{\begin{array}{llll}
h+2, & i\leq r,\\
h+1, &r<i\leq h,\\
r,& i=h +1,\\
0,& h+1<i.
\end{array}
\right.
$$

satisfies $\sum_i^t \lambda_i (T)=\sum_i^t  d^*_i (T) = 2m(T)$ and this is the maximum attainable over all threshold
graphs $T$ with $m(T)=m(G)$ edges. 
\end{enumerate}	
	Define Laplacian energy of graph as follows
	$$
	LE(G)=\sum_{i=1}^n \biggl|\mu_i (L(G))-\frac{2m(G)}{n}\biggr| .
	$$

	The main result of the paper~\cite{4} is the following
	\begin{theorem}
		For each spectrally threshold dominated graph G there exists a threshold graph with the same number of nodes and edges whose Laplacian energy is at least as large as that of G.
	\end{theorem}

	\section{Preliminary Remarks}
	
	Here we gather preliminary results that will be useful later.

	Let $\bar{G}= {[n]\choose 2}-G$ denote the complement of $G$. Then, (\cite{2}):
	$$
	\mu_i (L(G))=n-\mu_{n-i}(L(\bar{G})),\ i=1,\ldots ,n-1,
	$$

	The following duality result will be key in our work. It follows directly from the proof of Theorem 6 in \cite{1}, by including  $\Delta$'s with proper indices  in the calculation. 
	
	\begin{theorem}[\cite{3}]  \label{thm:dual}  For every graph $G$,
	$$
			\Delta_t(G)=\Delta_{n-t-1}(\bar{G})
			%			S_t(G)-{t+1\choose 2 }-m(G) = S_{n-t-1}(\bar G)-{n-t\choose 2}-m(\bar G).
		$$
		In particular, $G$ satisfies $BC_t$  if and only if $\bar{G}$ satisfies  $BC_{n-t-1}$.
	\end{theorem}
	
	On the other hand, once $G$ satisfy $BC_t$, the graph obtained by adding an isolated vertex, $G\cup\{v\}$,  trivially satisfy $BC_t$. Then, from Theorem \ref{thm:dual}, we conclude that the graph $G' = \overline{\bar{G}\bigcup\{ v\}}$ obtained by adding a dominating vertex $v$  satisfies $BC_{t+1}$:
	\begin{equation}\label{eq:threshold}
		\Delta_{t+1}(G') =\Delta_{t+1}(\overline{\bar{G}\cup\{v\}})=\Delta_{n-t-1}(\bar{G}\cup \{v\})=\Delta_{n-t-1}(\bar{G})=\Delta_t(G).   
	\end{equation}
	Given $G\subset {[n]\choose 2}$, we define the \textit{threshold family of $G$}, $\cal T(G)$, as the family of all graphs obtained from $G$ by adding complete or empty vertices. Note that the family of threshold graphs defined in the Introduction coincides with $\cal T(\emptyset)$. From Theorem~\ref{thm:dual} and ~equality (\ref{eq:threshold}) we conclude that $G$ satisfy Conjecture iff an element in $\cal T(G)$ does so. From this fact it follows
	
	%	give the first steps to our induction argument, whose proof we conclude in the end of the section:
	
	\begin{lemma}\label{lem:induction}
		Brouwer's Conjecture is valid for every $n$ and $t$ provided  that $BC_{t'}$ holds for every graph $G$ with $n'$ vertices where $t'={\frac{n'-1}{2}}$  if $n'$  is odd  or $t'$ equal to either ${\frac{n'-2}{2}}$ or ${\frac{n'}{2}}$ if $n'$ is even. 
	\end{lemma}
	
	We call the explicit $t'$s in Lemma \ref{lem:induction} as the \textit{middle $t$'s}. 	
	In what follows we will consider an inductive approach on $n$ to prove that $BC_t$ holds for the middle $t$'s, whenever it holds for middle $t$'s for graphs with fewer vertices.  To this end,  we remove one vertex of $G$  and   derive a special basis of $R^n$ where explicit bounds can be inferred. Recall the following formula for $L(G)$:
\begin{equation}
		(L(G)v,v)=\frac{1}{2}\sum_{(p,q)\in E} (v_p-v_q)^2. 
		\label{eq:L}
	\end{equation}	We have~(\cite{der},\ Cor 4.3.18)
	\begin{eqnarray*}&&
		S_t(G)=\max\left\{ \sum_{i=1}^t(L(G)x_i,x_i)\biggl| x_1,...,x_t, (x_i,x_j)=\delta_{ij}  \right\},\\ &&=\max\left(\hbox{tr}(L(G)_V)| V\hbox{ is a } t 
		\hbox{ dimensional subspace of }\ R^n\right\} \label{eq:S}	\\
		&&= \sum_{i=1}^t(L(G)z_i,z_i) 
	\end{eqnarray*}
	for $\{z_1,...,z_n\}$ an orthonormal set of eigenvectors corresponding to non-increasing eigenvalues of $L(G)$, and $z_n =z = (1/\sqrt{n},...,1/\sqrt{n})$.  
	
	 From the last equality we conclude that 
	\[S_t(G) = \sum_{i=1}^t(L(G)x_i,x_i)   \]
	for any orthonormal basis $\{x_1,...,x_t\}$ of $\hbox{span}\{z_1,...,z_t\}$.
	
	We have 
	\begin{eqnarray}&&
	\label{rr}
	S_t(G) = \max_{\{h_j,\ j\in [t]\}\in\hbox{ort} (n,t)} \sum_{i=1}^t(L(G)x_i,x_i) \leq   \max_{\{h_j,\ j\in [t]\}\in\hbox{ort} (n,t)}\sum_{i=1}^t(D h_i,h_i) + \max_{\{h_j,\ j\in [t]\}\in\hbox{ort} (n,t)}\sum_{i=1}^t(Ah_i,h_i)\\
	&& \leq\sum_{i=1}^{t}d_i + \sqrt{t\sum_{i=1}^n \alpha_i^2}\leq\sum_{i=1}^{t}d_i + \sqrt{nm}\leq \sum_{i=1}^{t}d_i + n\sqrt{n}, \nonumber
	\end{eqnarray}
	where $\alpha_i ,\ i\in [n]$ are eigenvalues of $A$ and $\hbox{ort}(n,t)$ is the family of sets of $t$ orthonormal vectors in $R^n$.
	
	\begin{lemma}\label{lem:xt}
		There exists an orthonormal basis $\{x_1,...,x_t\}$ of $\hbox{span}\{z_1,...,z_t\}$ and orthonormal basis 
		$\{x_{t+1},\ldots ,x_{n-1}\}$ \\ of $\hbox{span}\{z_{t+1},\ldots , z_{n-1}\}$ for any $t\in[n-1]$ 
		such that $x_i =(0,\ldots ,0,x_{i,i},x_{i,i+1},\ldots ,x_{i,n}),\ i\in [t],$\\
		$  x_{i} =(0,\ldots ,0,x_{i-t,i},x_{i-t-1,i},\ldots ,x_{n,i}),\ i\in [t+1,n-1]$.  	\end{lemma}
		We skip the usual proof of this Lemma, it contains the statement that one can choose the basis of such form in arbitrary subspace of dimension $t$ and $n-t-1$.    
			 
	From now we  fix a basis as in Lemma \ref{lem:xt} and denote it by $\{x_1,...,x_t,x_{t+1},...,x_{n-1}\}$. It is easy to see that $0\leq x_{1,1},x_{t+1,1}\leq\sqrt{\frac{n-1}{n}}$. 
This is because $\sum_{i=1}^n x_{i,j}^2 =1$ and $x_{n,j}=z_j =\frac{1}{\sqrt{n}}.$  We further assume $0<x_{1,1},x_{t+1,1}<\sqrt{\frac{n-1}{n}}$, since the extremal cases are easily dealt with. 	
	
  The existence of $x_1$ also allows our induction step. Let $x_1,...,x_t$ be as in Lemma \ref{lem:xt}. Given $G\subset {[n]\choose 2}$, consider the subgraph $G-\{1\}$ obtained by removing the first vertex of $G$, together with its edges. 
  
  We have 
	\begin{eqnarray*}&&
		S_t(G)=\sum_{i=1}^t(x_i,L(G)x_i)=\sum_{i=2}^{t} (x_i , L(G-\{1\})x_i)+\sum_{p; (1,p)\in E}\sum_{i=2}^{t}x_{i,p}^2+(x_1,L(G)x_1)\\ &&
		\leq S_{t-1}(G-\{1\})+\omega_1 +(x_1,L(G)x_1),
	\end{eqnarray*}
	where 
	\begin{equation}
		\omega_1 =\sum_{q:(1,q)\in E} \sum_{i=2}^{t} (x_{i,1}-x_{i,q})^2=\sum_{q:(1,q)\in E}\sum_{i=2}^t x_{i,q}^2 \leq d_1 .
		\end{equation}
	In particular, if $G-\{1\}$ satisfies $BC_{t-1}$, then $G$ satisfies $BC_t$ if
	\begin{equation}\label{eq:InductionStep}
		\omega_1 +(x_1,L(G)x_1)\leq t+d_1.
	\end{equation}
Equivalently, we can work with the complement graph, $\bar G$, and show that $BC_{\bar{t}}$ holds if $\bar{G}-\{1\}$ satisfies $BC_{\bar{t}-1}$ and
	\begin{equation}\label{eq:InductionStepDual}
		\bar\omega_1 +(x_{t+1},L(\bar{G})x_{t+1})\leq \bar  t+\bar d_1 .
	\end{equation}
	Here we take $x_{t+1}$ as the only vector with (possibly) non-zero first coordinate, and 
	\begin{equation}
		\bar{t}=n-1-t,\ \bar{d}_1= n-1-d_1,\  \ \bar{\omega}_1 =\sum_{q:(1,q)\in \bar{E}} \sum_{i=t+2}^{n-1} (x_{i,1}-x_{i,q})^2 = \sum_{q:(1,q)\in \bar{E}}\sum_{i=t+2}^{n-1} x_{i,q}^2 \leq \bar{d}_1 .
		\end{equation}
		
			It is easy to see that we can choose arbitrary $p\in [n]$ instead of the first coordinate in above consideration with substitution $1 \leftrightarrow p$ in the formulas, we use this consideration below several times.

	The key elements in the paper are the following bounds on $(L(G)x_1,x_1)$ 
	
	\begin{proposition}\label{prop:bound}
		Let $x_1$ be as in Lemma \ref{lem:xt} and $x_{1,1}> 0$. Then,
		\begin{eqnarray}
			\label{1q}
			(x_1 ,L(G)x_1)
			\leq\left\{\begin{array}{ll}
				d_1+\sqrt{d_1\frac{1-x_{1,1}^2}{x_{1,1}^2}}, & x_{1,1}^2 \geq \frac{d_1}{d_1 +1};\\
				\frac{nd_1}{n-1} +\sqrt{\frac{d_1\bar d_1}{n-1}\frac{1-\frac{n}{n-1}x_{1,1}^2}{x_{1,1}^2}},\ & x_{1,1}^2 <\frac{d_1}{d_1 +1}.
			\end{array}
			\right.
		\end{eqnarray}
		Likewise,
		\begin{eqnarray}
			\label{2q}
			(x_{t+1} ,L(\bar{G})x_{t+1}) \nonumber
			\leq\left\{\begin{array}{ll}
				\bar{d}_1+\sqrt{\bar{d}_1\frac{1-x_{t+1,1}^2}{x_{t+1,1}^2}}, & x_{t+1,1}^2 \geq \frac{\bar{d}_1}{\bar{d}_1 +1};\\
				\frac{n\bar{d}_1}{n-1} +\sqrt{\bar{d}_1\left(1-\frac{\bar{d}_1}{n-1}\right)\frac{1-\frac{n}{n-1}x_{t+1,1}^2}{x_{t+1,1}^2}},\ & x_{t+1,1}^2 <\frac{\bar{d}_1}{\bar{d}_1 +1}.
			\end{array}\right.\nonumber
		\end{eqnarray}
	\end{proposition}
	\begin{proof}
		By eventually replacing $x_1$ by $-x_1$, we assume that $x_{1,1}>0$.  
		
		We have
				\begin{eqnarray*}&&
			(x_1,L(G)x_1)x_{1,1}=d_1x_{1,1}-\sum_{p:(1,p)\in E}x_{t,p}\leq d_1x_{1,1}+\Biggl|\sum_{p:(1,p)\in E}x_{1,p}\Biggr|
		\end{eqnarray*}
		and
		$$
		\Biggl|\sum_{p:(1,p)\in E}x_{1,p}\Biggr|=\Biggl|\sum_{p:(1,p)\in \bar{E}}x_{1,p} +x_{1,1}\Biggr|\leq x_{1,1}+\Biggl|\sum_{p:(1,p)\in \bar{E}}x_{1,p}\Biggr|.
		$$
		
		Using Jensen inequality we obtain:
		$$
		\Biggl|\sum_{p:(1,p)\in E}x_{1,p}\Biggr|\leq\sqrt{d_1 x},\qquad \Biggl|\sum_{p:(1,p)\in \bar{E}}x_{1,p}\Biggr|\leq \sqrt{\bar{d}_1 (1-x_{1,1}^2 -x)},
		$$
		where $x=\sum_{p:(1,p)\in E}x_{1,p}^2$. 
		
		Therefore,
		\begin{eqnarray*}&&
		\Biggl|\sum_{p:(1,p)\in E}x_{1,p}\Biggr|\leq  
			\max_{x\in [0,1-x_{1,1}^2]} \min\left\{\sqrt{d_1 x},x_{1,1}+\sqrt{\bar{d}_1(1-x_{1,1}^2 -x)}\right\}\\ &&
			= \left\{\begin{array}{cc} \sqrt{d_1 (1-x_{1,1}^2)}, &~  x_{1,1}^2 \geq \frac{d_1}{d_1 +1};\\
				\frac{x_{1,1}d_1}{n-1} +\sqrt{\frac{d\bar d_{1}}{n-1}\left(1-\frac{n}{n-1}x_{1,1}^2\right)}, &  \hbox{otherwise}.
			\end{array}\right.
		\end{eqnarray*}
				
		The  first bound on $x_{1,1}^2$ is equivalent to 
		\[\sqrt{d_1 (1-x^2_{1,1})} \leq x_{1,1}, \]
		making $\sqrt{d_1 (1-x^2_{1,1})}$  the solution to the $\max\min$ problem. Otherwise,  since $x\mapsto \sqrt{\bar{d}_1(1-x^2_{1,1}-x)}$ is decreasing, 
		the $\max\min$ is achieved when
		\begin{equation}\label{eq:Le2equality}
			\sqrt{d_1 x}=x_{1,1}+\sqrt{\bar{d}_1(1-x_{1,1}^2 -x)}.	
		\end{equation}
		We manipulate this equation as follows:
		\begin{eqnarray*}&&
			(\sqrt{d_1 x}-x_{1,1})^2=\bar d_1(1-x_{1,1}^2 -x)\Leftrightarrow\\&&
			(n-1)x -2\sqrt{d_1}x_{1,1}\sqrt{x}+x^2_{t,1}-\bar d_1(1-x^2_{1,1})=0\Leftrightarrow\\&&
			\sqrt x = \frac{\sqrt{d_1}x_{1,1}}{n-1}+\sqrt{\left(\frac{\sqrt{d_1}x_{1,1}}{n-1}\right)^2-\frac{x^2_{1,1}-\bar d_1(1-x^2_{1,1})}{n-1} }\\&&
			=\frac{\sqrt{d_1}x_{1,1}}{n-1}+\sqrt{\frac{{d_1}x_{1,1}^2-  (n-1)x^2_{1,1}+(n-1)\bar d_1(1-x^2_{1,1})} {(n-1)^2} }\\&&
			=\frac{\sqrt{d_1}x_{1,1}}{n-1}+
			\sqrt{\frac{\bar d_1(1-x^2_{1,1}-\frac{1}{n-1}x_{1,1}^2)} {n-1} }\\&&
			=
			\frac{\sqrt{d_1}x_{1,1}}{n-1}+
			\sqrt{\frac{\bar d_1} {n-1}\left(1-\frac{n}{n-1}x^2_{1,1}\right) }.
		\end{eqnarray*}
		The result is concluded by multiplying the  last expression by $\sqrt{d_1}$.
	\end{proof}

	Before proceeding, we remark the following inequality that follows from the last proof.

	\begin{lemma}\label{cor:inn}
		Suppose $x_{1,1}^2<\frac{d_1}{d_1+1}$. Then, $\biggl|\sum_{p:(1,p)\in E}x_{1,p}\biggr|\leq<\sqrt{\frac{d_1\bar d_1}{n-1}}$.
	\end{lemma}
	\begin{proof}
		In the proof of Proposition \ref{prop:bound}, we concluded that
		\[\biggl|\sum_{p:(1,p)\in E}x_{1,p}\biggr|\leq		\frac{x_{1,1}d_1}{n-1} +\sqrt{\frac{d\bar d_{1}}{n-1}\left(1-\frac{n}{n-1}x_{1,1}^2\right)}.
		\]
		Proof follows from the observation that r.h.s. of last inequality is decreasing function of $x_{1,1}$ and hence achieved its maximum for $x_{1,1}>0$ at $x_{1,1}=0$.
			\end{proof}

	An extra inequalities are also needed. Recall that $x_1,x_{t+1}$ are the only vectors in $\{x_1,...,x_{n-1}\}$ with non-zero first coordinates. To motivate the next inequality, we also recall that the first vertex is complete if and only if the vector
	$$z=\Big(\sqrt{\frac{n-1}{n}}, -\frac{1}{\sqrt{n(n-1)}},...,-\frac{1}{\sqrt{n(n-1)}}   \Big) $$
	is in the span of $\{x_1,...,x_{t}\}$.
	
	 Next, we measure how much this vector does not belong to this $t$-subspace. 
	
	There exists $0<\lambda<1$ and  a vector  $y = (0, y_2\ldots ,y_n),\  \sum_{p=2}^n y_p =0,\  \sum_{p=2}^n y_p^2 =1$ such that
	\begin{eqnarray*}
		&& x_1= z\sqrt{\lambda } +  \sqrt{1-\lambda}y,\\
		&&  x_{t+1}= z\sqrt{1-\lambda} -\sqrt{\lambda}y.
	\end{eqnarray*}

		Further denote:
		\begin{equation}
			B=\frac{1}{2}(y,L(G)y)=\sum_{p<q,\ (p,q)\in E}(y_p -y_q)^2,\qquad \quad \bar{B}=\frac{1}{2}(y,L(\bar{G})y)=\sum_{p<q,\ (p,q)\in \bar{E}}(y_p -y_q)^2.
			\end{equation}	
		
		Then, using inequalities
		$$
		\biggl|\sum_{p:(1,p)\in E} y_{p}\biggr| \leq\max_{x\in [0,1]}\min \left\{\sqrt{d_1x},\sqrt{\bar{d}_1(1-x)}\right\}
		$$
		we have
		\begin{eqnarray}&&\nonumber
			(x_1,L(G)x_1) = \lambda d_1\frac{n}{n-1} + (1-\lambda )B -2\sqrt{ \lambda (1-\lambda )\frac{n}{n-1}}\sum_{p:(1,p)\in E} y_{p} \\
			&& \leq  \lambda d_1\frac{n}{n-1} +(1-\lambda )B +2\sqrt{\frac{n}{n-1} \lambda (1-\lambda )d_1 \left(1-\frac{d_1}{n-1}\right)}; \label{eq:tomax}\\ && \nonumber
			(x_{t+1},L(\bar{G})x_{t+1})= (1-\lambda) \bar{d}_1\frac{n}{n-1}+\lambda \bar{B} -2\sqrt{\lambda (1-\lambda )\frac{n}{n-1}}\sum_{p:(1,p)\in\bar{E}} y_{p}\\ \nonumber
			&& \leq  (1-  \lambda) \bar{d}_1\frac{n}{n-1} +\lambda \bar{B} +2\sqrt{\frac{n}{n-1} \lambda (1-\lambda )\bar{d}_1 \left(1-\frac{\bar{d}_1}{n-1}\right)}.
		\end{eqnarray}

		Optimization over $\lambda$ deliver the following bounds
		
		\begin{proposition} \label{prop:newinequality}
			Let $x_t$ be as above. Then,
			\begin{eqnarray*}&&
				(x_1,L(G)x_1)  \leq  \frac{d_1\frac{n}{n-1} +B}{2} +\frac{1}{2}\sqrt{\left(d_1\frac{n}{n-1} -B\right)^2 +4\frac{n}{n-1}d_1 \left(1-\frac{d_1}{n-1}\right)};\\ 
				&&(x_{t+1},L(\bar{G})x_{t+1})\leq  \frac{\bar{d}_1\frac{n}{n-1} +\bar{B}}{2} +\frac{1}{2}\sqrt{\left(\bar{d}_1\frac{n}{n-1} -\bar{B}\right)^2 +4\frac{n}{n-1}\bar{d}_1 \left(1-\frac{\bar{d}_1}{n-1}\right)}.
			\end{eqnarray*}
		\end{proposition}
		\begin{proof}
			We maximize the expression in \ref{eq:tomax} for $0< \lambda <1$. To this aim, we analyze the derivative of the expression with respect to $\lambda$:
			\begin{equation}\label{eq:equaltozero}
				d_1\frac{n}{n-1}-B+\sqrt{\frac{d_1\bar d_1}{n-1}}\frac{1-2\lambda}{\sqrt{\lambda(1-\lambda)}}.
			\end{equation}	
			Observe that the derivative goes to $+\infty$ and $-\infty$ as $\lambda$ goes to $0$ and $1$, respectively. 
			
			Therefore, we conclude that the maximum is in the interior.
			On the other hand setting expression \ref{eq:equaltozero} to zero gives:
			$$
			\lambda^2-\lambda + \frac{1}{4+A^2}=0, \qquad \quad A = \frac{\left(d_1\frac{n}{n-1}-B\right)(n-1)}{\sqrt{d_1\bar d_1 n}}.$$
			The maximum is achieved at:
			$$
			\lambda_{\pm} =\frac{1}{2}\left( 1\pm \frac{A}{\sqrt{4 +A^2}}\right).
			$$
			The proof is concluded by replacing $\lambda$ by $\lambda_\pm$  in \ref{eq:tomax}, observing that $\lambda_{\pm}=1-\lambda_{\mp}$.
		\end{proof}
		
		Using Proposition \ref{prop:newinequality}  and conditions~(\ref{eq:InductionStep}),  we conclude that if graph $G-\{1\}$ satisfies $\hbox{BC}_{t-1}$ and
		\begin{eqnarray}&&\label{bb1}
			\frac{d_1\frac{n}{n-1} +B}{2} +\frac{1}{2}\sqrt{\left(d_1\frac{n}{n-1} -B\right)^2 +4d_1\left(1-\frac{d_1}{n-1}\right)\frac{n}{n-1}} +\omega_1\leq t+d_1
			\end{eqnarray}
			then graph $G$ satisfies $BC_t$. Simular using condition~(\ref{eq:InductionStepDual}) 
			\begin{eqnarray}\label{bb2}
			&& \frac{\bar{d}_1\frac{n}{n-1} +\bar{B}}{2} +\frac{1}{2}\sqrt{\left(\bar{d}_1\frac{n}{n-1} -\bar{B}\right)^2 +4\bar{d}_1\left(1-\frac{\bar{d}_1}{n-1}\right)\frac{n}{n-1}}+\bar\omega_1  \leq \bar{t}+\bar{d}_1
		\end{eqnarray}	
		we conclude that if graph $\bar{G}-\{1\}$ satisfies $\hbox{BC}_{\bar{t}-1}$, then graph $\bar{G}$ satisfies $BC_{\bar t}$
		
				\section{Proof of Conjecture}
				
				We describe the key steps in the proof.
\begin{enumerate}
\item
First case.
Using (9) and assuming condition $x_{1,1}^2 \geq\frac{d_1}{d_1 +1}$ or $x_{t+1,1}^2 \geq\frac{\bar{d}_1}{\bar{d}_1 +1}$ we make inductive proof  $\hbox{BC}_t$ for graph $G$ or $\bar{G}$ assuming that $\hbox{BC}_{t-1}$ for $\hbox{G}-\{1\}$ is valid.    
\item
Second case.
We assume that there exists $p\in [n]$ s.t. $\omega_p \leq t(1-\delta )$. Because we can make permutation $p\leftrightarrow 1$ vertices of graph in arbitrary way, w.l.o.g. we set $p=1$. Assume also \begin{equation}
\label{fd}
x_{1,1}^2 \geq\frac{2}{n\delta^2},\ 1/5 >\delta >2n^{-1/3}
\end{equation}
and $t-B>\frac{2}{\delta}$.
At first we use inequality (6) which we reduce to (16) for one step inductive proof  $\hbox{BC}_t$ for graph $G$ under the condition that $\hbox{BC}_{t-1}$ for graph $G-\{1\}$ is true.
Next we consider the case $B\geq t-\frac{2}{\delta}$. Then we come to contradiction to the condition~(\ref{fd}). 
\item
Third case $\omega_q >t(1-\delta),\ q\in [n]$.

When
$$
m(\bar{G})\geq {t\choose 2}(1+ 3\delta ),
$$
we prove $\hbox{BC}_{t-1}$ for graph $\bar{G}$ directly by using bound (4).

In the case
$$
m(\bar{G})< {t\choose 2}(1+ 3\delta ),
$$
we first assume that there exist $p \in [n]$ s.t. $t(1+\delta )\geq d_p, \bar{d}_p \geq t(1-\delta )$. Next we show that there exist set ${\cal R}\subset [n]$ 
s.t. $\bar{d}_q \leq 7n\delta^{1/4},\ q\in {\cal R},$ 
$$
a =|{\cal R}|=\biggl[n(1-8\delta^{1/4})\left(1-\sqrt{1-\frac{1}{2(1-8\delta^{1/4})^2}}\right)\biggr].
$$
By permutation of vertices of graph w.l.o.g. we can assume that ${\cal R}=[a]$.
 We make $a$ steps of induction adding step by step $[a]$ vertices to the graph $G-[a]$ 
and for each step $i\in [a]$ we prove the $\hbox{BC}_{t-i}$ for the graph $G-[i]$ under the assumption that $BC_{t-i-1}$ is valid for graph $G-[i+1],\ i\in [a]$.
Choice of $a$ in (26) allows to prove $\hbox{BC}_{t-a}$ directly by using bound (4).  
\item
In the last case assumption is $d_q > t(1+\delta )$ or $\bar{d}_q \leq (1-\delta),\ q\in [n],\ \hbox{BC}_{t-1}$ for $\bar{G}$ is proved directly, using bound (4).  
\end{enumerate}

Next we use this scheme to demonstrate the proof in details.		
		
		We use induction on $n$ to prove BC and assume that BC is true for $n\leq 10^{24}$.  One can significantly improve this bound for $n$ by following the proof in this article more carefully.

		Let $\{x_i,\ i\in  [n-1]\}$ be the set of eigenvectors of $L(G)$. Considering Grassmannian frame 
		$F$ with row set $\{ x_i,\ i\in [t]\}$ and complement frame $\bar{F}$ with row set $\{x_i,\ i\in [t+1,n-1]\}$. 
				
		Note, that 
		$$
		x_{1,1}^2 +x_{t+1,1}^2 =\frac{n-1}{n}.
		$$
		
		W.l.o.g. we can assume that $x_{1,1}\in \left(0,\sqrt{(n-1)/n}\right)$.

		As a first step, we observe that the case $x_{1,1}^2\geq \frac{d_1}{d_1+1},\ \left( \hbox{respectively,}\ x_{t+1,1}^2\geq \frac{\bar{d}_1}{\bar{d}_1+1}\right)$   is easily discarded.
		\begin{lemma}
			Suppose that either $x_{1,1}^2\geq  \frac{d_1}{d_1+1}$ or $x_{t+1,1}^2\geq  \frac{\bar{d}_1}{\bar{d}_1+1}$. Then, BC holds for $G$.
		\end{lemma}	
		\begin{proof}
			To prove BC for $n$ and $x_{1,1}^2 \geq\frac{d_1}{d_1+1}$, assuming that it is true for $n-1$, it is sufficient to prove the inequality
			$$
			d_1 +\sqrt{d_1 \frac{1-x_{1,1}^2}{x_{1,1}^2}}+\omega_1 \leq d_1 +t
			$$
			or
			$$
			\omega_1 \leq t-\sqrt{d_1 \frac{1-x_{1,1}^2}{x_{1,1}^2}}.
			$$
			Last inequality is trivial, since 
			\[\sqrt{d_1 \frac{1-x_{1,1}^2}{x_{1,1}^2}}\leq \sqrt{d_1\frac{1-\frac{d_1}{d_1+1}}{\frac{d_1}{d_1+1}}}\leq 1.\] 
			The same consideration proves BC when $x_{t+1,1}^2 \geq\frac{\bar{d}_1}{\bar{d}_1+1}.$
		\end{proof}
%
%		A key point in the remaining of the proof is the symmetry of the square-rooted term with respect to taking duals:
%	\begin{gather}
%		\sqrt{\frac{d_1\bar d_1}{n-1}\frac{1-\frac{n}{n-1}x_{1,1}^2}{x_{1,1}^2}}=
%		\sqrt{\frac{n}{n-1}\frac{d_1\bar d_1}{n-1}\frac{\frac{n-1}{n}-x_{1,1}^2}{x_{1,1}^2}}=
%		\sqrt{\frac{n}{n-1}\frac{d_1\bar d_1}{n-1}\frac{x_{t+1,1}^2}{x_{1,1}^2}}.
%	\end{gather}
	
	Taking into account conditions  \ref{eq:InductionStep}, \ref{eq:InductionStepDual} and Proposition \ref{prop:bound} together we conclude that  $BC_t$ holds for $G$ if one of the following inequalities is true:
	\begin{eqnarray}
	&&
	\label{zz3}
	\frac{d_1}{n-1}+\sqrt{d_1\left(1-\frac{d_1}{n-1}\right)\frac{1-\frac{n}{n-1}x_{1,1}^2}{x_{1,1}^2}}+\omega_1 \leq t;\\
	\label{zz4}
	&&\frac{\bar{d}_1}{n-1}+\sqrt{\bar d_1\left(1-\frac{\bar d_1}{n-1}\right)\frac{1-\frac{n}{n-1}x_{t+1,1}^2}{x_{t+1,1}^2}}+\bar\omega_1 \leq \bar t.
\end{eqnarray}
For the remaining of the paper, we consider  $t=\frac{n}{2}$ when $n=2t$ and $t=\frac{n-1}{2}$ when $n=2t+1$.  Assume at first that $\omega_1 <t(1-\delta )$. 

Then using~(\ref{zz3}) we obtain the inequality
$$
{d_1\left(1-\frac{d_1}{n-1}\right)\frac{1-\frac{n}{n-1}x_{1,1}^2}{x_{1,1}^2}} \leq (t\delta -1)^2
$$
or
$$
x_{1,1}^2\geq\frac{s\bar{s}(n-1)}{s\bar{s}n+ (t\delta -1)^2},\  s=\frac{d_1}{n-1},\  \bar{s}=1-s.$$
The last inequality is satisfied if
\begin{equation}
\label{rr1}
x_{1,1}^2\geq \frac{2}{n\delta^2},\ \frac{1}{5}>\delta>3n^{-1/3}.
\end{equation}
It is left to consider the reverse condition:
\begin{equation}
\label{ff1}
x_{1,1}^2 <\frac{2}{n\delta^2}.
\end{equation}
In this case, we have:
\begin{equation} 
\label{rr2}
d_1 <\omega_1 +\sum_{p:(1,p)\in E} (x_{1,p}-x_{1,1})^2 <d_1x_{1,1}^2  +2\sqrt{d_1}x_{1,1} +1+\omega_1<\frac{1}{\delta^2} +\frac{2}{\delta}+1 +\omega_1< \frac{2}{\delta^2}+\omega_1 \leq t\left(1-\frac{\delta}{4}\right),\  3n^{-1/3}<\delta <10^{-1}.
\end{equation}	
	
	      Imposing condition~(\ref{bb1}) and using inequality $\omega_1 \leq d_1$, we obtain stronger condition
		 \begin{equation}
		 \label{ew}
		 \frac{d_1\frac{n}{n-1} +B}{2}+\frac{1}{2}\sqrt{\left(d_1 \frac{n}{n-1} - B\right)^2 +4d_1\left(1-\frac{d_1}{n-1}\right)\frac{n}{n-1}} +d_1\leq d_1 +t
		 \end{equation}
		 to prove that graph $G$ satisfies $\hbox{BC}_t$ if $G-\{1\}$ satisfies $\hbox{BC}_{t-1}$. 
		 
		 Hence 
		 $$
		 \sqrt{\left(d_1 \frac{n}{n-1} - B\right)^2 +4d_1\left(1-\frac{d_1}{n-1}\right)\frac{n}{n-1}} \leq 2t -d_1\frac{n}{n-1} -B.
		 $$
		 Assume that $ B\leq t -\frac{2}{\delta}$, then
		 $$
		 4d_1 \left( 1-\frac{d_1}{n-1}\right)\frac{n}{n-1} \leq \left(2t-d_1\frac{n}{n-1}-B\right)^2 - \left(d_1 \frac{n}{n-1} - B\right)^2.
		 $$
                  Hence we need to prove inequality
		 $$
		 d_1\left(1-\frac{d_1}{n-1}\right)\frac{n}{n-1} \leq \left(t-d_1 \frac{n}{n-1}\right) (t-B).
		 $$
		 
		 To satisfy last inequality it is sufficient to impose condition
		 $$
		 d_1 <t\left(1-\frac{\delta}{4}\right)
		 $$
            when $\delta >3n^{-1/3},\ t-B\geq\frac{2}{\delta}.$ 
		 
		 At last, if $B\geq t-\frac{2}{\delta},$ then $\bar{B}<n-t+\frac{2}{\delta}=\frac{n}{2} +\frac{2}{\delta},\ \bar{d}_1 > n-1 -t\left(1-\frac{\delta}{4}\right)\geq\frac{n-2}{2}\left(1+\frac{\delta}{4} \right).$ 
		 
		  From other side, 
		 \begin{eqnarray}\label{f2}&&
		 \bar{B}\geq \sum_{q:(1,q)\in\bar{E}} (x_{t+1,1}-x_{t+1,p})^2 \geq \bar{d}_1 x_{t+1,1}^2 -2|x_{t+1,1}|\sqrt{\bar{d}_1(1-x_{t+1,1}^2)}\\ \nonumber
		 &&>\frac{n-2}{2}\left(1+\frac{\delta}{4} \right)\left(\frac{n-1}{n}-\frac{2}{n\delta^2}\right)-2\sqrt{\left(\frac{2}{n\delta^2}+\frac{1}{n}\right)\frac{n}{2}\left(1+\frac{\delta}{4} \right)}\\
		 && \nonumber 
		 \geq\frac{n}{2}+n\frac{\delta}{8}-\left(2+\frac{\delta}{4}\right)\left(1+\frac{1}{\delta^2}\right) - \frac{2}{\delta}\left(1+\frac{\delta}{8} \right)\\
	 &&> \frac{n}{2}+\frac{2}{\delta},\ \hbox{where}\ \frac{1}{5}>\delta >3n^{-1/3}. \nonumber
		 \end{eqnarray}
		 This contradiction complete the proof in the case that  exists $p\in [n],$ s.t. $\omega_p \leq t(1-\delta)$. Here in the third inequality we use relation 
		 $$
		  x_{t+1,1}^2 =\frac{n-1}{n}-x_{1,1}^2$$ and in the forth inequality in~(\ref{f2}) we use the inequality
		 $$
		  2\sqrt{\left(\frac{2}{n\delta^2}+\frac{1}{n}\right)\frac{n}{2}\left(1+\frac{\delta}{4} \right)}< \frac{2}{\delta}\left(1+\frac{\delta}{8} \right).
		 $$  
		 
		 Next we assume that $\omega_q >t(1-\delta),\ q\in [n]$. Then $d_q \geq \omega_q  >t(1-\delta )$ and $\bar{d}_q \leq t(1+\delta),\ q\in [n]$. $\hbox{BC}_{\bar{t}},\ \bar{t}= t,t-1$ to be true for complement graph $\bar{G}$
		 it is sufficient to impose inequality
		 $$
		 \sum_{q=t+1}^{n-1}(x_{q},L(\bar{G})x_{q})) = \sum_{q=t+1}^{n-1}\mu_q (L(\bar{G}))\leq \sum_{q=t+1}^{n-1}\bar{d}_q +n\sqrt{n}\leq t^2(1+\delta) +n\sqrt{n}\leq m(\bar{G})+{\bar{t}\choose 2}.
		 $$
		 Here we use bound~(\ref{rr}).
		 
		 From the last inequality it follows that $\hbox{BC}_{\bar{t}}$ is true for $\bar{G}$ if the number of edges $m(\bar{G})$ satisfies the inequality
		 $$
		 m(\bar{G})\geq{t\choose 2}(1+3\delta).
		 $$
		 			
		 	 Assume now that there exists $p\in[n]$ s.t. $t(1+\delta )> d_p,\ \bar{d}_p \geq t(1-\delta )$.  
			 
		 Taking into account that $\omega_p >t(1-\delta )$ we have 
		 $$
		 \sum_{q: (p,q)\in \bar{E}}\sum_{i=1}^t x_{i,q}^2 <t-t(1-\delta)=t\delta.
		 $$
		 Assuming uniform distribution on the set $\{ q: (p,q)\in \bar{E}\}$,
		 $$ E\left(\sum_{i=1}^t x_{i,q}^2\right) \leq  \frac{t\delta}{\bar{d}_p} \leq \frac{t\delta}{t(1-\delta)}=\frac{\delta}{1-\delta}.
		 $$
		 Using Markov inequality $P(X> CE(X))\leq\frac{1}{C},\  C>0$ and choosing $C=\delta^{-1/2}$,
		 we have
		 $$
		 P\left(\sum_{i=1}^t x_{i,q}^2 >\frac{\sqrt{\delta}}{1-\delta }\right)<\sqrt{\delta}.
		 $$
		 Hence
		 $$
		 \sum_{i=1}^t x_{i,q}^2 <\frac{\sqrt{\delta}}{1-\delta}
		 $$
		 where $q\in J$ for some set $J\subset \{ q: (p,q)\in \bar{E}\},\ |J|>(1-\sqrt{\delta})(1-\delta )t.$  Note also that $\omega_q \leq t$.
		 
		 Hence for the arbitrary $q\in J$
		 \begin{eqnarray*}&&
		 d_q \leq \sum_{r: (r,q)\in E} \sum_{i=1}^t (x_{i,r} -x_{i,q})^2 \leq \sum_{r: (r,q)\in E} \sum_{i=1}^t x_{i,q}^2 +\sum_{r: (r,q)\in E} \sum_{i=1}^t x_{i,r}^2
		 -2\sum_{r: (r,q)\in E} \sum_{i=1}^t x_{i,q}x_{i,r}\\ && \leq t +d_q \frac{\sqrt{\delta}}{1-\delta }+2\sum_{r: (r,q)\in E} \sqrt{\sum_{i=1}^t x_{i,q}^2}\sqrt{\sum_{i=1}^t x_{i,r}^2} \leq t+ d_q \frac{\sqrt{\delta}}{1-\delta }+2d_q \sqrt{\frac{\sqrt{\delta}}{1-\delta }}
		 \end{eqnarray*}
		 or
		 \begin{equation}
		 \label{rr3}
		 d_q \leq \frac{t}{1-\frac{\sqrt{\delta}}{1-\delta } -2\sqrt{\frac{\sqrt{\delta}}{1-\delta }}} \leq t(1+3\delta^{1/4}) ,\ q\in J,\ \delta >3n^{-1/3}.
		 \end{equation}
		
		 Because $d_p>\omega_p >t(1-\delta)> t(1-4\delta^{1/4}),\  p\in [n]$, we have
		 inequalities
		 $$
		 t(1-4\delta^{1/4})<d_p,\bar{d}_q < t(1+4\delta^{1/4}),\ \omega_q >t(1-\delta),\ q\in J.
		 $$
		 Thus
		 $$
		\sum_{q\in [n]\setminus J} \bar{d}_ q \leq    2m(\bar{G}) - (1-4\delta^{1/4}) t^2(1-\delta )(1-\sqrt{\delta}) \leq 2{t\choose 2}(1+3\delta) - (1-4\delta^{1/4}) t^2(1-\delta )(1-\sqrt{\delta})< 3n^2 \sqrt{\delta}.
		$$
		W.l.o.g. we can assume $|J|=[t(1-\delta)(1-\sqrt{\delta})]$.
		Assuming uniform distribution on the set $[n]\setminus J$, we have
		 $E(\bar{d}_q )=\frac{\sum_{q\in  J} \bar{d}_ q}{|[n]\setminus J|} <\frac{3n^2\sqrt{\delta}}{n-t(1-\delta)(1-\sqrt{\delta)}} <7n\sqrt{\delta}$. Using Markov inequality we have
		$$
		P\left( \bar{d}_q \geq 7n\delta^{1/4} \right)< \delta^{1/4}.
		$$
                 Thus there exists set $I\subset [n]\setminus J$, $|I| >n- t(1-\delta^{1/4})(1-\sqrt{\delta})(1-\delta)>t(1+\delta^{1/4})$ s.t. $\bar{d}_q < 7n\delta^{1/4},\ q\in I$ and hence
                 $d_q \geq n(1-7\delta^{1/4}),\ q\in I.$ 
                 W.l.o.g. we can assume that  $I=\lfloor t(1+\delta^{1/4})\rfloor$ and it is sufficient here to assume that $7\delta^{1/4}\leq 1/10$.  
                 
                 Last inequality and inequality $3n^{-1/3} <\delta$, which is imposed in~(\ref{rr1}),~(\ref{rr2}),~(\ref{f2}),~(\ref{rr3})  leads to the condition  $3n^{-1/3}< \delta  < (70)^{-4}$. Hence $3n^{-1/3}  < (70)^{-4}$. The last inequality to be true it is sufficient to impose condition $n>10^{24}$.    
                 \bigskip
                 
       We consider  $a$ coordinates ${\cal R}$ from the set $I$. W.l.o.g. we can assume that ${\cal R}=[a]$, where
       $$a=\biggl[n(1-8\delta^{1/4})\left(1-\sqrt{1-\frac{1}{2(1-8\delta^{1/4})^2}}\right)\biggr].$$  Justification of the choice  of  $a$ we make later. When passing step $i\in [a-1]$ we renumber vertices  of graph $G-[i-1]$ as follows $i\to i-1$ skipping first $i-1$ positions in graphs $G-[i-1]$ and $\bar{G}-[i-1]$. On this way we redefine $t(i)=t-i$ vectors $x_{t+1},\ldots ,x_{n-1}$
       as follows $x_{t+i}=(0,\dots ,0, x_{t+i,i},\ldots ,x_{t+i,n})\to  \tilde{x}_{t+i} =( x_{t+i,1}, \ldots , x_{t+i, n-i+1}).$  Complement set of orthonormal vectors of length $n-i+1$ we denote
       $\tilde{x}_1 ,\ldots ,\tilde{x}_t $.
       \bigskip
       
       $\bigcirc$\bigskip \ \ {\bf\hbox{Starting point for the process, described below and implemented} "a" \hbox{times}.} \bigskip
       
       Using inequality~(\ref{zz4}) on $i$-th step and taking into account the inequalities $\bar{\omega}_1 <\bar{d}_1 \leq 7n\delta^{1/4}$ and $t(i) = t -i$ we obtain the inequality
       \begin{eqnarray}
       \label{f1}&&
       \bar{d}_1 \left( 1-\frac{\bar{d}_1}{n-i}\right) \frac{1-\frac{n-i+1}{n-i}\tilde{x}^2_{t+i,1}}{\tilde{x}_{t+i,1}^2} \leq\left(t-1-i-7n\delta^{1/4}\right)^2 
       \end{eqnarray}
       Because $i\leq a =\biggl[\frac{3n}{8(1-8\delta^{1/4})}\biggr]$, we relax bound~(\ref{f1}) to
       $$
       \tilde{x}_{t+i,1}^2 \geq\frac{7n\delta^{1/4}}{\left(\frac{1}{14}n-7n\delta^{1/4}\right)^2} >\frac{7\cdot 14^2 \delta^{1/4}}{n}.
       $$
       Assume now, that the opposite inequality is valid $\tilde{x}_{t+i,1}^2 <\frac{7\cdot 14^2 \delta^{1/4}}{n}$. Then we repeat considerations  starting from equation~(\ref{ew}) for $\bar{d}_1(G-[i]) \to \bar{d}_1 (G-[i-1]) \leq 7n\delta^{1/4} ,\ t\to t-i,\ \bar{B}(G-[i])\to \bar{B}(G-[i-1])$.
       
       According the induction process, we impose the inequality~(\ref{bb2})
       \begin{eqnarray}\label{t5}&&
       \frac{\bar{d}_1 \frac{n-i+1}{n-i} +\bar{B}}{2}+\frac{1}{2}\sqrt{\left(\bar{d}_1\frac{n-i+1}{n-i}-\bar{B}\right)^2 +4\bar{d}_1 \left(1-\frac{\bar{d}_1}{n-i}\right)\frac{n-i+1}{n-i}}+\bar{\omega}_1  \\
       && \leq \bar{d}_1 +t-i.\nonumber
       \end{eqnarray}
       Assume that $\bar{B}\leq t-i-1$.
       Making transformations of the last inequality, and using inequality $\bar{\omega}_1 \leq \bar{d}_1$, we impose stronger inequality
       $$
       \bar{d}_1 \left( 1-\frac{\bar{d}_1}{n-i}\right)\frac{n-i+1}{n-i} \leq (t-i-\bar{B})\left( t-i-\bar{d}_1 \frac{n-i+1}{n-i}\right).
       $$
       We strength last inequality and obtain the bound
       $$
       8n\delta^{1/4}\leq(t-i-\bar{B})(t-i-7n\delta^{1/4}) 
       $$
        or
        $$
        \bar{B} \leq t -i -\frac{8n\delta^{1/4}}{t-i-7n\delta^{1/4}}. 
        $$
        Assume now the opposite inequality
        $$
\bar{B} >t -i -\frac{8n\delta^{1/4}}{t-i-7n\delta^{1/4}}
        $$ 
        or
        \begin{equation}
        \label{vc}
    B\leq t+1+\frac{8n\delta^{1/4}}{t-i-7n\delta^{1/4}}.
    \end{equation}
        From other side we have $d_1 >n(1-7\delta^{1/4})$ and $\tilde{x}^2_{1,1} >\frac{n-i}{n-i+1}-\frac{7\cdot 14^2\delta^{1/4}}{n}>
        1-\frac{2}{n},\ \delta <(70)^{-4}$.
        
        At the end we show that  inequality~(\ref{vc}) could not be satisfied, and we come to contradiction:
        \begin{eqnarray*}
        && B\geq \sum_{p:(1,p)\in E(G-[i-1])} (\tilde{x}_{1,1}-\tilde{x}_{1,p})^2 >d_1 \tilde{x}_{1,1}^2 -2|\tilde{x}_{1,1}|\sqrt{d_1} >n(1-7\delta^{1/4}) \left(1-\frac{2}{n}\right)\\
        && -2\sqrt{n(1-7\delta^{1/4}) \left(1-\frac{2}{n}\right)} >\frac{3}{2}t.
        \end{eqnarray*}
        Last inequality contradict to inequality~(\ref{vc}) for $n>10^{24}$.

                 Next we make above proof procedure  (from sign $\bigcirc$ ) for graph $\bar{G}-[i-1]$ and set of orthonormal vectors $x_{t+i}, \ldots , x_{n-1}$. Step by step deleting vertex $i\in I$ from $\bar{G}$ on $i$-th step and assuming by induction

                 that $\hbox{BC}_{t-1-i}$ is true for graph $\bar{G}- [i]$ and set of vectors $x_{t+i+1}, \ldots , x_{n-1}$ of length $n-i=2t-i$ 
                 and as before proving that $\hbox{BC}_{t-i}$
                 is true for graph $\bar{G}-[i-1]$ and set
                 of vectors  $x_{t+i}, \ldots , x_{n-1}$.

                 We make $a$ steps of this induction process and obtain from graph $\bar{G}-[a]$ and set of vectors  $x_{t+a+1}, \ldots ,x_{n-1}$ of length $ n-a=2t-a$, graph $\bar{G}$ and set of vectors $x_{t+1}, \ldots ,x_{n-1}$ of length $n$. The
                 complement to graph $\bar{G}-[a]$ is graph $G-[a]$ and complement set of orthonormal vectors $x_1 ,\ldots , x_t$ of length $n-a$. The $\hbox{BC}_{t-a-1}$ is true for $\bar{G}-[a]$ iff   $\hbox{BC}_{t}$ is true for graph $G-[a]$.
                 
                  Remind that we assume that                 
                 $$
                 m(\bar{G}) \leq {t\choose 2}(1+3\delta) 
                 $$
                 or
                 $$
                {2t+1\choose 2}-{t\choose 2}\leq m(G) \geq 3{t\choose 2}(1-\delta) .
                 $$
                
                 $\hbox{BC}_{t}$ for $G-[a]$ is obviously true if $m(\bar{G}-[a])\leq {t+1\choose2}$. Because $d_q >n-1-\bar{d}_q \geq n-1-7n\delta^{1/4}>n(1-8\delta^{1/4})$, we have
                 $$
                 m(G-[a])\leq m(G)-na(1-8\delta^{1/4}\leq  {2t+1\choose 2}-{t\choose 2} - na(1-8\delta^{1/4} )+{a\choose 2}\leq {t+1\choose 2}.
                 $$
                 The last condition to be true  it is sufficient to impose condition
                 \begin{equation}
                 \label{ew2}
                 a=\biggl[n(1-8\delta^{1/4})\left(1-\sqrt{1-\frac{1}{2(1-8\delta^{1/4})^2}}\right)\biggr].
                 \end{equation}
                 Note that $|I|=[t(1+\delta^{1/4})] >a$.                 
                 
                  Assume now that $d_q > t(1+\delta ),\ i\in[n]$. Then $\bar{d}_q \leq t(1-\delta )$.
		$\hbox{BC}_{t-1}$ is true if
		$$
		\sum_{q=1}^{t-1}\mu_i (L(\bar{G})) \leq \sum_{i=1}^{t-1} \bar{d}_q+\sqrt{2m(\bar{G})t} \leq t(t-1)(1-\delta ) +\sqrt{nm(\bar{G})}  \leq t(t-1)(1-\delta ) +n\sqrt{n}\leq m(\bar{G})+\frac{n(n-2)}{8},
		$$
		which is true when $m(\bar{G})> {t\choose 2}$, otherwise $\hbox{BC}_{t-1}$ is trivially true for graph $\bar{G}$.  
                This completes the proof. \bigskip

		{\bf Acknowledgement}
\bigskip

Research supported by Extrato de contrato N114/2019, N20/2021, Processo N23089.101560/2018-91.  The article was
written while the first author was visiting Unifesp, Brazil. He would like to
thank the Department of Mathematics and also colleagues for their kind atmosphere to execute  this work.

\end{document}